\documentclass[10pt]{amsart}

\usepackage{amssymb,amsmath,amsthm}
\usepackage{graphicx,a4wide}

\theoremstyle{plain}

\usepackage{graphicx,tikz}
\newtheorem{theorem}{Theorem}

\newtheorem*{proposition}{Proposition}

\newtheorem*{lemma}{Lemma}

\theoremstyle{definition}

\theoremstyle{remark}

\newcommand{\vol}{\operatorname{vol}}

\begin{document}

\title[]{A Metric Sturm-Liouville Theory in two dimensions}

\author[]{Stefan Steinerberger}
\address{Department of Mathematics, Yale University}
\email{stefan.steinerberger@yale.edu}
\thanks{This work is supported by the NSF (DMS-1763179) and the Alfred P. Sloan Foundation.}

\begin{abstract}
A central result of Sturm-Liouville theory (also called the Sturm-Hurwitz Theorem) states that if $\phi_k$ is a sequence of eigenfunctions of a second order
differential operator on the interval $I \subset \mathbb{R}$, then any linear combination satisfies a uniform bound on the roots
$$ \# \left\{x \in I:\sum_{k \geq n}{ a_k \phi_k(x)} = 0 \right\} \geq n-1.$$
We provide a sharp (up to logarithmic factors) generalization to two dimensions: let $(M,g)$ be a compact two-dimensional manifold (with or without boundary), let
$(\phi_k)$ denote the sequence of eigenfunctions of a uniformly elliptic operator $-\mbox{div}(a(\cdot) \nabla)$ (with Dirichlet or Neumann boundary conditions). Then, 
for any linear combination of eigenfunctions above a certain index $n$,
$$ f = \sum_{k \geq n}{a_k \phi_k} ~ \mbox{we have} \quad \mathcal{H}^1 \left\{ x: f(x) = 0\right\}  \gtrsim_{} \frac{\sqrt{n}}{\sqrt{\log{n}}}     \log \left(n \frac{\|f\|_{L^2(M)}}{\|f\|_{L^1(M)}}  \right)^{-1/2}    \frac{\|f\|_{L^1(M)}}{\| f \|_{L^{\infty}(M)}} .$$
Examples on $M=\mathbb{T}^2$ and $M=\mathbb{S}^2$ shows that this is optimal up to the logarithmic factors. The proof is using optimal transport and a new inequality for the Wasserstein metric $W_p$: if
$f(x)dx$ and $g(x)dx$ are two absolutely continuous measures on a two-dimensional domain $M$ with continuous densities and the same total mass, then, for all $1 \leq p <\infty$,
$$ W_p(f(x)dx, g(x) dx) \cdot \mathcal{H}^1 \left\{x \in M: f(x) = g(x) \right\} \gtrsim_{M,p} \frac{\|f-g\|_{L^1(M)}^{1+1/p}}{\|f-g\|_{L^{\infty}(M)}}.$$
\end{abstract}

\maketitle

\vspace{0pt}

\section{Introduction}
\subsection{Sturm-Liouville Theory}
Sturm-Liouville theory dates back to seminal papers from 1836  \cite{liouville, sturm, sturm2} and is concerned with oscillation properties of eigenfunctions of operators
$$ H = -\frac{d}{dx}\left(a(x) \frac{d}{dx} \right) + b(x) \qquad \mbox{on an interval}~(a,b)$$
where $a(x), b(x) > 0$ are bounded away from 0 (this is not an exhaustive description of Sturm-Liouville theory, we refer to Galaktionov \& Harwin \cite{gal} or Zettl \cite{anton}). Sturm proved that there exists a discrete set of parameter $(\lambda_n)_{n=1}^{\infty}$ (the eigenvalues of the Sturm-Liouville
operator $H$) and an associated sequence of solutions $(\phi_n)_{n=1}^{\infty}$ that form an orthogonal basis in $L^2(a,b)$
and that the number of their roots is completely rigid (called the Sturm Oscillation Theorem in textbooks).

\begin{quote}
\textbf{Weak Sturm Oscillation Theorem.} $\phi_n$ has $n-1$ roots in $(a,b)$.
\end{quote} 

 However, both Sturm and Liouville originally proved a \textit{much} stronger result (Sturm being the first to establish the result, Liouville then gave a different proof).
 That stronger result is not very well known (we could not find it any textbook, for example) and
reads as follows.

\begin{quote}
\textbf{Original Sturm Oscillation Theorem.} For any integers $m \leq n$ and any set of coefficients $a_m, a_{m+1}, \dots, a_n$ such that not all of them are 0, the function
$$ \sum_{k=m}^{n}{a_k \phi_k} \qquad \mbox{has at least}~m-1~\mbox{and at most}~n-1~\mbox{roots in}~(a,b).$$
\end{quote} 
This theorem seems to have been largely forgotten:  B\'erard \& Helffer \cite{berard} in a beautiful recent paper chronicle the decay of knowledge (as well as describing the original
proofs in modern language). Lord Rayleigh still called it 'a beautiful theorem' in 1877 but the book of Courant \& Hilbert already does not refer to it all (the suspicion in \cite{berard} being that Courant and Hilbert
did not consult the original papers but instead relied on a 1917 book of B\^{o}cher \cite{bocher} that also does not mention it).
The special case $\phi_n(x) = \sin{nx}$ is sometimes known as the Sturm-Hurwitz theorem after being stated by Hurwitz \cite{hurwitz} in 1903 (who explicitly refers to Sturm). This result has a particularly beautiful physical proof due to Polya \cite{pol}. A quantitative version was given by \cite{stein3}. An analogue of the Sturm-Hurwitz Theorem for the Fourier transform on the real line $\mathbb{R}$ was conjectured by Logan \cite{logan} and proved by Eremenko \& Novikov \cite{eremenko, erem2}.

\subsection{Sturm-Liouville Theory in higher dimensions.}
Let $M$ be a compact manifold (with or without boundary), let $H= - \mbox{div}(a(x)\cdot \nabla)$ be a uniformly elliptic operator of Laplacian-type and let $\phi_n$ denote the eigenfunctions of $H$ forming a basis of $L^2(M)$ with either Dirichlet or Neumann boundary conditions. It is difficult to speak of a Sturm-Liouville theory in higher dimension because there is more than one way of interpreting the one-dimensional statement.  Broadly speaking, the existing lines of research fall into two branches.\\

\begin{enumerate}
\item \textbf{Topological Sturm-Liouville Theory} interprets the one-dimensional results as a statement about the number of connected domains after we remove the zero set, i.e.
$$ \mbox{the number of connected components of} \qquad M \setminus \left\{x: \phi_n(x) = 0 \right\}.$$
Topological Sturm-Liouville theory has been of substantial interest to a number of people, including Courant (whose student Herrmann gave a flawed argument in that direction \cite{herrmann}), Gelfand (as recalled by Arnold \cite{arn3}) and Arnold \cite{arn, arn2}. It was Arnold who discovered that a straightforward generalization of the Sturm Oscillation theorem to the sphere would contradict results surrounding the topology of algebraic plane curves related to Hilbert's 16th problem. Gelfand proposed an approach that, as was known to himself, was unfortunately restricted to one dimension (this proof was recently reconstructed and completed by B\'{e}rard \& Helffer \cite{berard2}). The topological investigation of a single eigenfunction $\phi_k$ on a two-dimensional domain was initiated by Courant \cite{courant} who proved that $\phi_n$ has at most $n$ nodal domains. This was later improved by Pleijel \cite{plei} to $0.7n$ for $n$ sufficiently large, see also \cite{bour, stein0}. It is conjectured \cite{ber, polt} that the sharp constant might be $2/\pi \sim 0.63$ which would be attained on $\mathbb{T}^2$. Some recent results in that direction are \cite{berh, jayjay, lena}.
Topological Sturm-Liouville theory seems wide open.

\begin{center}
\begin{figure}[h!]
\begin{tikzpicture}[scale = 1.6]
\draw[ultra thick, draw=black] (0,0) rectangle ++(2,1);
\draw[ultra thick, draw=black] (2.5,0) rectangle ++(2,1);
\draw[ultra thick, draw=black] (0,1.5) rectangle ++(2,1);
\draw[ultra thick, draw=black] (2.5,1.5) rectangle ++(2,1);
\node at (-0.5, 0.5) {$n=2$};
\node at (-0.5, 2) {$n=1$};
\node at (1, 3) {Single eigenfunction};
\node at (3.5, 3) {Orthogonal to};
\node at (3.5, 2.7) {low frequencies};
\node at (0.9, 2.2) {Weak Sturm};
\node at (1.05, 1.9) {Oscillation Thm.};
\node at (2.5+0.9, 2.2) {Strong Sturm};
\node at (2.5+1.05, 1.9) {Oscillation Thm.};
\node at (1, 0.7) {$\mathcal{H}^1(\phi_n = 0) \gtrsim \sqrt{n}$};
\node at (1, 0.2) {Br\"uning (1978, \cite{brun})};
\node at (2.5+1.05, 0.5) {this paper};
\end{tikzpicture}
\caption{Metric Sturm-Liouville theory in one and two dimensions.}
\end{figure}
\end{center}

\item \textbf{Metric Sturm-Liouville Theory} interprets the one-dimensional result as a statement about the size of the $(n-1)-$dimensional Hausdorff measure of the zero set. A famous conjecture of S.-T. Yau \cite{yau} states that if $-\Delta \phi_n = \lambda \phi_n$, then
$$  \mathcal{H}^{n-1} \left\{x: \phi_n(x) = 0 \right\} \sim \sqrt{\lambda}.$$
This has attracted considerably amount of interest with contributions by Br\"uning \cite{brun} (who first established the lower bound in $n=2$ dimensions), Chanillo \& Muckenhoupt \cite{chan}, Colding \& Minicozzi \cite{colding}, R. T. Dong \cite{dong} (who first established the upper bound $\lambda^{3/4}$ in two dimensions), Donnelly \cite{donn}, Donnelly \& Fefferman \cite{don, don2} (who proved the conjecture for analytic metrics), Q. Han \& F.-H. Lin \cite{lin0}, Q. Han, R. Hardt \& F.-H. Lin \cite{lin1}, Hardt \& Simon \cite{hardt} (who proved an upper bound in all dimensions), Hezari \& Wong \cite{hezari}, Hezari \& Sogge \cite{hezari2}, Jerison \& Lebeau \cite{jer}, F.-H. Lin \cite{lin2}, Logunov \& Malinnikova \cite{mal}, Mangoubi \cite{mangoubi}, Nadirashvili \cite{nadirashvili}, Sogge \& Zelditch \cite{sog1, sog2} and the author \cite{stein1}. The lower bound was recently established in all dimensions by Logunov \cite{logunov}.
\end{enumerate}

\section{Main Results}

\subsection{Metric Sturm-Liouville theory.}
We will now state our main result. Let $(M,g)$ be a compact two-dimensional manifold (with or without boundary) and let 
$$H = -\mbox{div}(a(x) \cdot \nabla)$$
be a uniformly elliptic second-order operator equipped with either Dirichlet or Neumann boundary conditions (in case of Neumann boundary conditions, we would ask that $\partial M$ has some degree of regularity). Let
$\phi_n$ denote the sequence of eigenfunctions of $H$. Weyl's theorem implies that the eigenvalues scale like $\lambda_n \sim_M n$, where the implicit constant depends only on the area of $M$. 
Our main result shows that any linear combination of eigenfunctions above a certain frequency has an unavoidable degree of vanishing: the sum of many oscillating functions is still oscillating.

\begin{theorem} If $f \in C^0(M)$ is orthogonal to the first $n$ eigenfunctions, i.e. if 
$ f = \sum_{k \geq n}{a_k \phi_k}$,
then we have the following estimate on the length of its nodal line
 $$\mathcal{H}^1 \left\{ x: f(x) = 0\right\}  \gtrsim_{M} \frac{\sqrt{n}}{\sqrt{\log{n}}}     \log \left(e+n \frac{\|f\|_{L^2(M)}}{\|f\|_{L^1(M)}}  \right)^{-1/2}    \frac{\|f\|_{L^1(M)}}{\| f \|_{L^{\infty}(M)}} .$$
\end{theorem}
Note that $f$ being orthogonal to the first $n$ eigenfunctions corresponds to $n$ equations being satisfied, the Theorem is a statement about an infinite-dimensional space.
The result has the optimal scaling in $n$ (up to a factor of $\log{n}$): consider $M = \mathbb{T}^2$ and the function 
$$ f(x,y) = \sin{(\sqrt{n} x)} \qquad \mbox{satisfying} \quad  \mathcal{H}^1 \left\{ x: f(x) = 0\right\} \sim \sqrt{n}, \quad \| f \|_{L^{p}(M)} \sim 1.$$
One could wonder whether better results are possible on other manifolds. We are able to obtain the
following restrictions on estimates of this type: if there is an estimate of the form
$$  \mathcal{H}^1 \left\{ x: f(x) = 0\right\}  \gtrsim_{(M,g)}  n^{\alpha} \left( \frac{\|f\|_{L^1(M)}}{\| f \|_{L^{\infty}(M)}}\right)^{\beta} \qquad \mbox{then}~\alpha \leq \frac{1}{2}~\mbox{and }\beta \geq \frac{1}{2}.$$ 
The first statement, $\alpha \leq 1/2$, is suggested by known results on the statistics of arithmetic random waves on the torus $\mathbb{T}^2$ (see e.g. \cite{aurich}). Taking random linear combinations
of eigenfunctions at eigenvalue $\sim \sqrt{n}$ produces functions whose nodal length statistics still obey Yau's heuristic $\sim \sqrt{n}$ but for which $\|f\|_{L^{\infty}(\mathbb{T}^2)} \leq (\log{n})^{\gamma} \|f\|_{L^1(\mathbb{T}^2)}$ with high likelihood. This suggests $\alpha \leq 1/2$. We will obtain both statements from the following construction that seems to be new.

\begin{proposition} For some universal $c > 0$, every $n \in \mathbb{N}$ and every $0 < t \leq (c n)^{-1}$, there exists a function $f \in C^2(\mathbb{S}^2)$ orthogonal to the first $n$ eigenfunctions of the Laplacian $-\Delta_{\mathbb{S}^2}$ satisfying (up to logarithmic factors)
$$ \mathcal{H}^1 \left\{ x: f(x) = 0\right\} \sim n \sqrt{t} \qquad \mbox{as well as} \qquad  \frac{\|f\|_{L^1(M)}}{\| f \|_{L^{\infty}(M)}} \sim n t.$$
\end{proposition}
We do not know whether the estimate in Theorem 1 holds for some $1/2 \leq \beta \leq 1$ and consider this an interesting problem. 
We also mention an earlier result of the author \cite{stein1} that established a metric Sturm-Liouville theory with suboptimal exponents in all dimensions (but weaker than Theorem 1 in dimension $n=2$) by a different method: that result in $d$ dimensions reads
$$ \mathcal{H}^{d-1} \left\{x:f(x) =0\right\} \gtrsim_{M}      \frac{ n^{1/d}}{(\log{n})^{d/2}}  \left( \frac{ \|f\|_{L^1}}{\|f\|_{L^{\infty}}} \right)^{2 - \frac{1}{d}}  .$$

\subsection{A Wasserstein inequality.} Our argument is based on a new geometric inequality for optimal transport that may be of independent interest. 
The Wasserstein metric is a notion of distance between measures introduced in the late 1960s \cite{dob, wasser} and is now a foundational concept in optimal transport, probability theory and partial differential equations  \cite{otto, villani}. We define the $p-$Wasserstein distance between two measures $\mu$ and $\nu$ on a domain $M$ via
$$ W_p(\mu, \nu) = \left( \inf_{\gamma \in \Gamma(\mu, \nu)} \int_{M \times M}{ |x-y|^p d \gamma(x,y)}\right)^{1/p},$$
where $| \cdot |$ is the distance and $\Gamma(\mu, \nu)$ denotes the collection of all measures on $M \times M$
with marginals $\mu$ and $\nu$, respectively (also called the set of all couplings of $\mu$ and $\nu$). 
The special case $p=1$ is particularly nice:
in many settings we have Monge-Kantorovich duality (see e.g. \cite{villani})
$$ W_1(\mu, \nu) = \sup\left\{ \int_{M}{ f d\mu} - \int_{M}{ f d\nu}: f~\mbox{is 1-Lipschitz} \right\}.$$
The 1-Wasserstein distance or Earth Mover's Distance is the total amount of work ($=\mbox{distance}\times \mbox{mass})$ required to
move $\mu$ to $\nu$. 
Let now $\Omega \subset \mathbb{R}^2$ be a bounded domain or a compact two-dimensional manifold and let $\mu, \nu$ be two measures on $\Omega$ satisfying $\mu(\Omega) = \nu(\Omega)$. 
 We assume that $\mu$ and $\nu$ are
absolutely continuous with respect to the Lebesgue measure and that their densities are continuous function and are thus given by $\mu = g(x) dx$ and $\nu = h(x) dx$ (with $g,h \in C^0(\Omega)$). We introduce 
the function
$$ f(x) = g(x) - h(x)$$
and show that it captures some information about the transportation cost $W^p(\mu, \nu)$: if it changes sign along a large one-dimensional set, then we would expect that the measures
are fairly well mixed and that it should be cheap to move one to the other.
 However, if the zero set $\left\{x:f(x) = 0\right\}$ is rather short, then it seems like it would be quite
expensive to move on the other.

\begin{center}
\begin{figure}[h!]
\begin{tikzpicture}[scale = 1.5, yscale=2]
\draw [ultra thick] (-1,0) -- (5.5,0);
\draw [thick] (0,0) to[out=70, in =120] (1,0);
\draw [thick] (1,0) to[out=70, in =120] (2,0);
\draw [thick] (2,0) to[out=70, in =120] (3,0);
\draw [thick] (3,0) to[out=70, in =120] (4,0);
\draw [dashed, thick] (0.5,0) to[out=70, in =120] (1.5,0);
\draw [thick, dashed] (1.5,0) to[out=70, in =120] (2.5,0);
\draw [thick, dashed] (2.5,0) to[out=70, in =120] (3.5,0);
\draw [thick, dashed] (3.5,0) to[out=70, in =120] (4.5,0);
\draw [thick] (5, 0.4) -- (5.4, 0.401); \node at (5.5, 0.4) {$f$};
\draw [thick] (5, 0.2) -- (5.4, 0.201); \node at (5.5, 0.2) {$g$};
\end{tikzpicture}
\caption{A toy picture of Theorem 2 in one dimension: $f(x)dx$ is fairly cheap to transport to $g(x)dx$ in $W^1$ but $f(x) - g(x)$ also changes sign often.}
\end{figure}
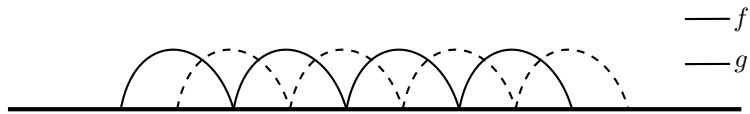
\end{center}

 We prove a geometric inequality making this notion precise: it has the flavor of an uncertainty principle, either the zero set is large or the transport is not cheap.

\begin{theorem} Let $f \not \equiv 0$. Then we have, for all $1 \leq p < \infty$,
$$ W_p(\mu, \nu) \cdot \mathcal{H}^1\left( \left\{x \in \Omega: f(x) = 0\right\} \right) \gtrsim_{p,\Omega} \frac{\|f\|_{L^1(\Omega)}^{1+ 1/p^{}}}{\|f\|_{L^{\infty}(\Omega)}}.$$
\end{theorem}

We are not aware of any such result existing in the literature with the exception of a one-dimensional version of the statement of a similar flavor that was established by the
author \cite{stein2}: for all continuous $f:\mathbb{T} \rightarrow \mathbb{R}_{}$  with mean value 0
$$ (\mbox{number of roots of}~f) \cdot \left( \sum_{k=1}^{\infty}{ \frac{ |\widehat{f}(k)|^2}{k^2}}\right)^{\frac{1}{2}}  \gtrsim \frac{\|f\|^{2}_{L^1(\mathbb{T})}}{\|f\|_{L^{\infty}(\mathbb{T})}}.$$
We note that the sum on the left-hand side is merely the Sobolev norm $H^{-1}$ which has connections to optimal transport via the infinitesimal expansion of the Wasserstein distance
$W^2$ \cite[\S 7.6]{villani} and an inequality of Peyr\'{e} \cite{pey}. Note added in print: Amir Sagiv and the author \cite{amir} have proved a variant of the result using a completely different idea that results in a similar estimate on $(0,1)^d$ where we show that, for all $1 \leq p < \infty$,
$$  W_p(\mu, \nu) \cdot \mathcal{H}^{d-1}\left\{x \in (0,1)^d: f(x) = 0 \right\} \gtrsim_{p,d}   \left( \frac{\|f\|_{L^1}}{\|f\|_{L^{\infty}}} \right)^{3 - \frac1d + \frac1p} \|f\|_{L^1}.$$
We note that this result is weaker than Theorem 2 for $d=2$.

\section{Proofs}

\subsection{A Geometric Lemma.} This section describes a simple geometric statement that is at the heart of the argument and also firmly restricts its applicability to two dimensions. The statement is so elementary that it is likely to be stated in the Literature in some form somewhere.
\begin{lemma}
Let $\Omega \subset \mathbb{R}^2$ be a connected domain. If $\varepsilon \leq |\Omega|^{1/2}/8$
then, for some universal $c>0$,
$$  \left| \left\{x \in \mathbb{R}^2 \setminus \Omega: d(x, \Omega) \leq \varepsilon\right\} \right| \leq c \cdot \varepsilon \cdot |\partial \Omega|.$$
\end{lemma}
We set $|\partial \Omega| = \infty$ if the boundary of $\Omega$ is not rectifiable, the statement is then trivially true.

\begin{proof} We first assume that $\Omega$ is simply connected. Let us define the set
$$A = \left\{x \in \mathbb{R}^2 \setminus \Omega: d(x, \Omega) \leq \varepsilon\right\}.$$
 We will define a transport of the Lebesgue measure of $A$ to the $\mathcal{H}^1$ measure on $\partial \Omega$ in the following way: pick a tiny square $Q$ in $A$ and distribute its Lebesgue measure evenly over $ \partial \Omega \cap \left\{x \in \mathbb{R}^2: d(x,Q) \leq 2 \varepsilon\right\}.$

\begin{center}
\begin{figure}[h!]
\begin{tikzpicture}[scale = 1]
\filldraw (0,0) circle (0.04cm);
\node at (0.2, -0.2){$x$};
\draw [dashed, thick] (0,0) circle (1cm);
\draw [dashed, thick] (0,0) circle (2cm);
\draw [ultra thick] (-2,-2) to[out=30, in = 270] (-0.8, 0) to[out=90, in =0] (-2, 1);
\node at (-2.4, -1.9) {$\partial \Omega$};
\end{tikzpicture}
\caption{The boundary exits the $2\varepsilon$-circle and is at least somewhere at most distance $\varepsilon$ from $x$. The triangle inequality forces the length to be $\geq 2\varepsilon$.}
\end{figure}
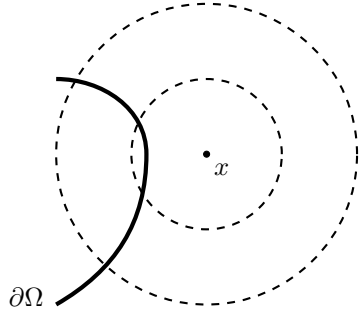
\end{center}

We claim that this leads to an even distribution (in the sense of the Radon-Nikodym derivative of the measure so created with respect to the Hausdorff measure $\mathcal{H}^1$ being bounded from above) because, uniformly for all $x \in A$,
$$ \mathcal{H}^1 \left( \partial \Omega \cap \left\{y \in \mathbb{R}^2: d(x,y) \leq 2 \varepsilon\right\} \right) \geq 2 \varepsilon.$$
This statement follows from the simple fact that $\varepsilon \leq |\Omega|^{1/2}/8$ implies that $\Omega$ is not contained in any $2 \varepsilon-$ball and the boundary thus has to leave the domain (see Fig. 3).
This implies that the Radon-Nikodym derivative of the induced measure w.r.t. to the Hausdorff measure is bounded from above by 
$\lesssim \varepsilon$ and this implies the result. It remains to discuss the case of multiply connected domains $\Omega$: in that case we argue in the very same way but transport only to the boundary of each connected component.
\end{proof}
As follows easily from the argument, the Lemma easily translates to general two-dimensional manifolds (the constant then, naturally, depends on the properties of the metric).
The argument has a somewhat vague similarity with Besicovitch's proof \cite{bes} of the systolic inequality in two dimensions. The Euclidean concentration inequality says that
among all sets $\Omega \subset \mathbb{R}^n$ with fixed volume, the ball has the smallest $\varepsilon-$enlargement $ \left| \left\{x \in \mathbb{R}^n \setminus \Omega: d(x, \Omega) \leq \varepsilon\right\} \right|$. It could be interesting to understand whether the ball is also extremal (in the other direction) for fixed surface area: is it true that among all sets with fixed boundary size the ball has the largest $\varepsilon-$enlargement? For $\varepsilon \rightarrow 0^+$ this question presumably reduces to a known statement about curvature.

\subsection{Proof of Theorem 2}

\begin{proof} 
We first give the proof for $p=1$ since that is the only case relevant in our application of Theorem 2 to Theorem 1 and then detail the necessary modifications for $p>1$.
We may assume w.l.o.g. that $\mu = g(x) dx$ and $\nu = h(x) dx$ and we assume that both $g$ and $h$ are continuous. Let $D \subset \Omega$ denoted a
connected component of $\left\{x \in \Omega: g(x) > h(x) \right\}$. We do not know anything about the transport plan that moves $g$ to $h$ but any such transport
plan has to at least transport the superfluous measure outside of that connected component $D$.
We use $\delta$ to denote the amount of $L^1-$mass that has to be transported outside of $D$,
$$ \delta = \int_{D}{g(x) - h(x) dx} = \int_{D}{f(x) dx}.$$
The question is how much it has to be transported: if we just move it barely outside of $D$, then it is going to be on a big pile and will not solve the problem. The best case
is if there is a big deficiency outside the domain and that
$$ g(x) - h(x) = -\|f\|_{L^{\infty}} \qquad \mbox{just outside of}~D.$$
The natural scale $\varepsilon$ on which we have to transport a typical particle then necessarily satisfies
$$  \left| \left\{x \in \mathbb{R}^2 \setminus D: d(x, \Omega) \leq \varepsilon\right\} \right| \cdot  \|f\|_{L^{\infty}} \gtrsim \delta.$$
The geometric Lemma implies that either $\varepsilon \gtrsim |D|^{1/2}$ or
$$ \delta \lesssim   \left| \left\{x \in \mathbb{R}^2 \setminus D: d(x, D) \leq \varepsilon\right\} \right| \|f\|_{L^{\infty}} \lesssim \varepsilon \cdot |\partial D| \cdot \|f\|_{L^{\infty}}.$$
It turns out that we can assume to deal with the second case since that is the weaker one; indeed, the second case can never be as good as the first case since, using the isoperimetric inequality,
$$ \frac{\delta}{\|f\|_{L^{\infty}}} \frac{1}{|\partial D|} =  \frac{1}{\|f\|_{L^{\infty}}} \frac{1}{|\partial D|} \int_{D}{f(x) dx} \leq     \frac{|D|}{|\partial D|} \lesssim |D|^{1/2}$$
and thus the arising lower bound on $\varepsilon$ can never exceed $|D|^{1/2}$. We thus assume, henceforth, that
$$ \varepsilon \gtrsim \frac{\delta}{\|f\|_{L^{\infty}}| |\partial D|}.$$
This shows that the cost of transporting the mass exceeding expectations outside of $D$ has
$$ \mbox{a}~W^1~\mbox{cost of at least} \gtrsim \varepsilon \delta \gtrsim \frac{\delta^2}{|\partial D|} \frac{1}{\|f\|_{L^{\infty}}}.$$
Let us now assume that $\left\{x \in \Omega: g(x) > h(x)\right\}$ has $n$ connected components $D_1, \dots, D_n$ (the subsequent estimates will not depend on $n$), then
$$ W^1(\mu, \nu) \gtrsim \frac{1}{\|f\|_{L^{\infty}}}\sum_{k=1}^{n}{ \frac{\|f\|_{L^1(D_i)}^2}{|\partial D_i|} }.$$
We conclude the argument with an application of the Cauchy-Schwarz inequality: since $f$ has mean value 0, we have
\begin{align*}
\frac{\|f\|_{L^1(\Omega)}}{2} &=  \sum_{k=1}^{n}{\|f\|_{L^1(D_i)}} =  \sum_{k=1}^{n}{\frac{\|f\|_{L^1(D_i)}}{|\partial D_i|^{1/2}} |\partial D_i|^{1/2}} \\
&\leq   \left( \sum_{k=1}^{n}{\frac{\|f\|^2_{L^1(D_i)}}{|\partial D_i|^{}}}\right)^{1/2} \left(   \sum_{k=1}^{n}{ |\partial D_i|^{}} \right)^{1/2} \\
&=    \left( \sum_{k=1}^{n}{\frac{\|f\|^2_{L^1(D_i)}}{|\partial D_i|^{}}}\right)^{1/2} 
\left(\mathcal{H}^1 \left\{x \in \Omega: f(x) = 0\right\} \right)^{1/2}\\
&=    \| f\|_{L^{\infty}}^{1/2} \left( \frac{1}{\|f\|_{L^{\infty}}} \sum_{k=1}^{n}{\frac{\|f\|^2_{L^1(D_i)}}{|\partial D_i|^{}}}\right)^{1/2} 
\left(\mathcal{H}^1 \left\{x \in \Omega: f(x) = 0\right\} \right)^{1/2}\\
&\leq    \| f\|_{L^{\infty}}^{1/2}  W_1(\mu, \nu)^{1/2} 
\left(\mathcal{H}^1 \left\{x \in \Omega: f(x) = 0\right\} \right)^{1/2}
\end{align*}
and therefore
$$ W^1(\mu, \nu) \cdot \mathcal{H}^1 \left\{x \in \Omega: f(x) = 0\right\} \gtrsim \frac{\|f\|^2_{L^1(\Omega)} }{\|f\|_{L^{\infty}}}.$$
The relevant changes for $p > 1$ are minimal: the lower bound on the transport cost is 
$$ W^p_p(\mu, \nu) \gtrsim \frac{1}{\|f\|_{L^{\infty}}^p}\sum_{k=1}^{n}{ \frac{\|f\|_{L^1(D_i)}^{p+1}}{|\partial D_i|^p} }$$
and the Cauchy-Schwarz inequality can be replaced by H\"older's inequality
\begin{align*}
\frac{\|f\|_{L^1(\Omega)}}{2} &=  \sum_{k=1}^{n}{ \frac{ \|f\|_{L^1(D_i)}}{ |\partial D_i|^{\frac{p}{p+1}}} |\partial D_i|^{\frac{p}{p+1}}  }   \leq
\left(  \sum_{k=1}^{n}{ \frac{ \|f\|^{p+1}_{L^1(D_i)}}{ |\partial D_i|^{p}}} \right)^{\frac{1}{p+1}} \left( \sum_{k=1}^{n}{ |\partial D_i|} \right)^{\frac{p}{p+1}}\\
&= \left(  \sum_{k=1}^{n}{ \frac{ \|f\|^{p+1}_{L^1(D_i)}}{ |\partial D_i|^{p}}} \right)^{\frac{1}{p+1}} \left(\mathcal{H}^1 \left\{x \in \Omega: f(x) = 0\right\} \right)^{\frac{p}{p+1}}
\end{align*}
and thus
$$  \sum_{k=1}^{n}{ \frac{ \|f\|^{p+1}_{L^1(D_i)}}{ |\partial D_i|^{p}}} \gtrsim_p \frac{\|f\|_{L^1}^{p+1}}{  \left(\mathcal{H}^1 \left\{x \in \Omega: f(x) = 0\right\} \right)^{p}}.$$
Altogether, we obtain
$$ W^p_p(\mu, \nu) \gtrsim \frac{1}{\|f\|_{L^{\infty}}^p}  \sum_{k=1}^{n}{ \frac{ \|f\|^{p+1}_{L^1(D_i)}}{ |\partial D_i|^{p}}} \gtrsim \frac{\|f\|_{L^1}^{p+1}}{\|f\|_{L^{\infty}}^p}  \frac{1}{  \left(\mathcal{H}^1 \left\{x \in \Omega: f(x) = 0\right\} \right)^{p}}$$
and this is the desired result.
\end{proof}

\subsection{Proof of Theorem 1.}
\begin{proof} Let us fix $f$ as
$$f = \sum_{k \geq n}{a_k \phi_k}.$$
We will apply Theorem 2 with 
$$ g(x) = \max\left\{ f(x), 0 \right\} \qquad \mbox{and} \qquad h(x) = -\min\left\{f(x), 0\right\}.$$
The desired result will then follow from Theorem 2 with $p=1$ and showing that
$$ W_1(g(x)dx,h(x)dx) \lesssim \frac{ \sqrt{ \log\left(  n \frac{\|f\|_{L^2(M)}}{\|f\|_{L^1(M)}}   \right) }}{ \sqrt{n}} \|f\|_{L^1}.$$
The estimate on the Wasserstein distance has previously been obtained by the author \cite{stein2} at a slightly greater level of generality, we give
a streamlined argument for $n=2$ dimensions. We decompose the function with respect to eigenfunctions (and note that the Weyl asymptotic in two dimensions 
is simply $\lambda_n \sim n$ and thus, in what follows, $\lambda \sim n$)
$$ f = \sum_{\lambda_k \geq \lambda}{ \left\langle f, \phi_k \right\rangle \phi_k}$$
and note that the solution of the heat equation $(\partial_t - \mbox{div}(a(x) \cdot \nabla))f_t= 0$ with $f_0 = f$ as initial conditions is explicitly given via convolution with the heat kernel or, alternatively, by diagonalization with eigenfunctions,
$$ f_t(x) = \int_{M}{ p(t,x,y) f(y) dy}  = \sum_{\lambda_k \geq \lambda}{ e^{-\lambda_k t} \left\langle f, \phi_k\right\rangle \phi_k}$$
and, in particular, we can estimate its size in $L^1$ from above by
\begin{align*} \|f_t\|^2_{L^1(M)} &= \left\| \sum_{\lambda_k \geq \lambda}{ e^{-\lambda_k t} \left\langle f, \phi_k\right\rangle \phi_k} \right\|^2_{L^1(M)} \leq \vol(M) \left\| \sum_{\lambda_k \geq \lambda}{ e^{-\lambda_k t} \left\langle f, \phi_k\right\rangle \phi_k} \right\|^2_{L^2(M)} \\
&= \vol(M) \sum_{\lambda_k \geq \lambda}{ e^{-2 \lambda_k t} \left| \left\langle f, \phi_k\right\rangle \right|^2} \lesssim_{(M,g)} e^{-2\lambda t}  \sum_{\lambda_k \geq \lambda}{ \left| \left\langle f, \phi_k\right\rangle \right|^2}
\lesssim_{(M,g)} e^{-2 \lambda t} \|f\|_{L^2(M)}^2
\end{align*}
At the same time, we can interpret diffusion as a microscopic process that takes a particular measure $\delta_x$ and spreads it according to the probability distribution $p(t,x,\cdot)$. The transport cost of this process can thus be bounded by, appealing to classical Gaussian bounds  \cite{aronson} in two dimensions
$$ p(t,x,y) \leq \frac{c_1}{t^{}} \exp \left( -\frac{|x-y|^2}{c_2 t} \right),$$
from above by
\begin{align*}
 W^1( \delta_x, p(t,x,y)dx) &\leq \int_{M}{p(t,x,y)\|x-y\| dy}\\
&\lesssim \int_{M}{ \frac{|x-y|}{t}\exp \left( -\frac{|x-y|^2}{c_2 t} \right) dy} \lesssim_{M} \sqrt{t}.
\end{align*}
The second argument is a trivial estimate: if we have two functions $h_1, h_2$ having the same total mass, then the total transport cost satisfies
$$ W^1(h_1(x)dx, h_2(x) dx) \leq \mbox{diam}(M) \|h_1 - h_2\|_{L^{\infty}}$$
Combining all these estimates, we see that
$$ W_1(g(x)dx,h(x)dx) \lesssim \sqrt{t} \|f\|_{L^1(M)}  +  e^{- \lambda t} \|f\|_{L^2(M)}.$$
Setting 
$$t = \lambda^{-1} \log{\left( \frac{\lambda^{\frac{1}{2}} \|f\|_{L^2(M)}}{\|f\|_{L^1(M)}}\right)} \qquad \mbox{yields} \qquad  W_1(g(x)dx,h(x)dx \lesssim
\frac{ \sqrt{ \log\left(  \lambda \frac{\|f\|_{L^2(M)}}{\|f\|_{L^1(M)}}  \right) }   }{ \sqrt{\lambda}} \|f\|_{L^1}.$$
\end{proof}

There is a corresponding upper bound on the $W_p$ distance for $p>1$, however, in terms of bounds on the
nodal set all the exponents cancel in such a way that the arguent does not yield a different result for different values of $p$, the implicit constant decays as $p \rightarrow \infty$
and can thus not be used to remove the logarithmic factor.

\subsection{Proof of the Proposition.}
\begin{proof}
The construction makes use of a rather recent construction of Bondarenko, Radchenko \& Viazovska \cite{bond, bond2} (answering a long-standing question in spherical designs): they prove the existence of a set $n$ points $\left\{x_1, \dots, x_n\right\}$ on $\mathbb{S}^2$ that are $n^{-1/2}-$separated such that for all polynomials $p:\mathbb{R}^3 \rightarrow \mathbb{R}$ up to degree $\sim \sqrt{n}$ (a vector space with dimension $\sim n$)
$$ \frac{1}{n} \sum_{k=1}^{n}{p(x_k)} = \frac{1}{|\mathbb{S}^2|} \int_{\mathbb{S}^2}{p(x) dx}.$$
However, polynomials restricted to $\mathbb{S}^2$ are exactly linear combinations of eigenfunctions of the Laplacian which means that the signed measure
$$ \mu = - \frac{n}{|\mathbb{S}^2|} + \sum_{k=1}^{n}{\delta_{x_k}}$$
is orthogonal to all eigenfunctions of the Laplacian up to eigenvalue $\sim \sqrt{n}$ (which contains the first $n$ eigenfunctions). Clearly, this measure is not a continuous function. We apply the heat equation up to time $t$ and obtain
$$ f_t = e^{t\Delta} \mu =  - \frac{n}{|\mathbb{S}^2|} + \sum_{k=1}^{n}{e^{t\Delta} \delta_{x_k}}.$$
Since the heat equation is diagonalized by eigenfunctions of the Laplacian, the heat flow preserves all orthogonality properties and $f_t$ is also orthogonal to the first $n$ eigenfunctions. The classical short-time asymptotics for the heat equation now imply that 
$$ \left(e^{t \Delta} \delta_x\right)(y) \sim \frac{1}{4\pi t} \exp\left(- \frac{\|x-y\|^2}{4t} \right).$$
The $n^{-1/2}-$separation of the points allows us to run the heat equation up to time $t \sim n^{-1}$ before heat balls start to overlap substantially. The asymptotics then immediately imply that, ignoring logarithmic factors,
$$ \mathcal{H}^1\left( \left\{x \in \mathbb{S}^2: f_t(x) = 0\right\}\right) \sim n \sqrt{t}$$
while
$$ \| f_t\|_{L^1(\mathbb{S}^2)} \sim n \qquad \mbox{and} \qquad  \| f_t\|_{L^{\infty}(\mathbb{S}^2)} \sim \frac{1}{t}.$$
\end{proof}

The construction is not limited to $\mathbb{S}^2$ (because \cite{bond2} is done in full generality on $\mathbb{S}^d$). The same argument yields, for all $d \geq 2$, all $n \geq 2$ and for all $0< t \leq c_d n^{-2/d}$ a function $f_t \in C^{\infty}(\mathbb{S}^d)$ such that $f$ is orthogonal to the first $n$ eigenfunctions of $-\Delta_{\mathbb{S}^{d-1}}$ and
$$ \mathcal{H}^{d-1}\left( \left\{x \in \mathbb{S}^d: f_t(x) = 0\right\}\right) \sim n t^{\frac{d-1}{2}}$$
as well as
$$ \| f_t \|_{L^1(\mathbb{S}^d)} \sim n \qquad \mbox{and} \qquad  \| f_t \|_{L^{\infty}(\mathbb{S}^2)} \sim \frac{1}{t^{d/2}}.$$
This shows that any $d-$dimensional statement in metric Sturm-Liouville theory of the form
$$  \mathcal{H}^{d-1} \left\{ x: f(x) = 0\right\}  \gtrsim_{(M,g)}  n^{\alpha} \left( \frac{\|f\|_{L^1(M)}}{\| f \|_{L^{\infty}(M)}}\right)^{\beta}$$
has to satisfy
$$ \alpha \leq \frac{1}{d} \qquad \mbox{and} \qquad \beta \geq \frac{d-1}{d}.$$
This scaling, coupled with the Weyl asymptotic, suggests the natural conjecture
$$  \mathcal{H}^{d-1} \left\{ x: f(x) = 0\right\}  \gtrsim_{(M,g)}  \sqrt{\lambda_n} \left( \frac{\|f\|_{L^1(M)}}{\| f \|_{L^{\infty}(M)}}\right)^{1-\frac1d}.$$


\begin{thebibliography}{10}

\bibitem{aronson} D. Aronson,  Non-negative solutions of linear parabolic equations, Ann. Sci. Norm. Sup. 22 (1968), 607--694.

\bibitem{arn} V. I. Arnol'd, The branched covering $CP2 \rightarrow S^4$, hyperbolicity and projective topology. (Russian) 
Sibirsk. Mat. Zh. 29 (1988), no. 5, 36--47, 237; translation in 
Siberian Math. J. 29 (1988), no. 5, 717--726 (1989). 

\bibitem{arn2}  V.I. Arnold, Topological problems in wave propagation theory and topological economy
principle in algebraic geometry, Third Lecture by V. Arnold at the Meeting in
the Fields Institute Dedicated to His 60th Birthday, Fields Inst. Commun., 1997.

\bibitem{arn3} V. Arnold. Topological properties of eigenoscillations in mathematical physics. Proceedings
of the Steklov Institute of Mathematics 273 (2011) 25--34.

\bibitem{aurich} R. Aurich, A. B\"acker, R. Schubert and M. Taglieber,
Maximum norms of chaotic quantum eigenstates and random waves. 
Phys. D 129 (1999), no. 1-2, 1--14.

\bibitem{ber} P. B\'erard, In\'egalit\'es isop\'erim\'etriques et applications. Domaines nodaux des fonctions propres. Goulaouic--Meyer--Schwartz Seminar, 1981/1982, Exp. No. XI, 10 pp., Ecole Polytech., Palaiseau, 1982.

\bibitem{berard} P. B\'{e}rard and B. Helffer, Sturm's Theorem on Zeros of Linear Combinations of Eigenfunctions, arXiv: 1706.08247

\bibitem{berard2} P. B\'{e}rard and B. Helffer, Sturm’s theorem on the zeros of sums of
eigenfunctions: Gelfand’s strategy implemented, arXiv:1807.03990.

\bibitem{berh} P. B\'{e}rard and B. Helffer, Courant-Sharp Eigenvalues for the Equilateral Torus, and for the Equilateral Triangle, Letters in Mathematical Physics 106 (2016), p. 1729 -- 1789.

\bibitem{bes} A. S. Besicovitch, On two problems of Loewner. J. London Math. Soc. 27, (1952). 141--144. 

\bibitem{bocher} M. B\^{o}cher. Lecons sur les methodes de Sturm dans la theorie des equations
differentielles lineaires et leurs developpements modernes. Gauthier-Villars et
Cie, \'Editeurs. Paris 1917.

\bibitem{bond} A. Bondarenko, D. Radchenko, M. Viazovska, Optimal asymptotic bounds for spherical designs, Annals of Mathematics 178 (2), 443--452, (2013).
\bibitem{bond2} A. Bondarenko, D. Radchenko, M. Viazovska, Well-separated spherical designs.  Constr. Approx. 41 (2015), no. 1, 93--112. 

\bibitem{bour}  J. Bourgain, On Pleijel’s Nodal Domain Theorem, IMRN 13 (2013), 1--7.

\bibitem{brun} J. Br\"uning, \"Uber Knoten Eigenfunktionen des Laplace-Beltrami Operators. Math. Z. 158:15--21, 1978.

\bibitem{chan} S. Chanillo, B. Muckenhoupt,  Nodal geometry on Riemannian manifolds. J. Diff. Geom. 34:85--91. 

\bibitem{colding} T. H. Colding and W. P. Minicozzi II, Lower bounds for nodal sets of eigenfunctions. Comm. Math. Phys. 306:  777--784, 2011.

\bibitem{courant} R. Courant. Ein allgemeiner Satz zur Theorie der Eigenfunktionen selbstadjungierter Differentialausdr\"ucke,
Nachr. Ges. G\"ottingen (1923), 81---4.

\bibitem{dob} R. Dobrusin, 
Definition of a system of random variables by means of conditional distributions. 
Teor. Verojatnost. i Primenen. 15 (1970), 469--497. 

\bibitem{dong} R. T. Dong, Nodal sets of eigenfunctions on Riemann surfaces. J. Differential Geom. 36 (1992), no.
2, 493--506.

\bibitem{donn} H. Donnelly,
Nodal sets for sums of eigenfunctions on Riemannian manifolds. 
Proc. Amer. Math. Soc. 121 (1994), no. 3, 967--973. 

\bibitem{don} H. Donnelly and C. Fefferman, 
Nodal sets of eigenfunctions on Riemannian manifolds. 
Invent. Math. 93 (1988), no. 1, 161--183. 


\bibitem{don2} H, Donnelly and C. Fefferman, Nodal sets for eigenfunctions of the Laplacian on surfaces. J. Amer. Math. Soc. 3:332--353.

\bibitem{eremenko} A. Eremenko and D. Novikov,
Oscillation of Fourier integrals with a spectral gap. 
J. Math. Pures Appl. (9) 83 (2004), no. 3, 313--365. 

\bibitem{erem2} A. Eremenko and D. Novikov,
Oscillation of functions with a spectral gap. 
PNAS (2004), 5872--5873. 

\bibitem{gal} V. Galaktionov, P. Harwin, 
Sturm's theorems on zero sets in nonlinear parabolic equations. Sturm-Liouville theory, 173--199, Birkh\"auser, Basel, 2005. 



\bibitem{lin0} Q. Han, F.-H. Lin, 
On the geometric measure of nodal sets of solutions.
J. Partial Differential Equations 7 (1994), no. 2, 111--131. 

\bibitem{lin1} Q. Han, R. Hardt and F.-H. Lin, 
Geometric measure of singular sets of elliptic equations. 
Comm. Pure Appl. Math. 51 (1998), no. 11-12, 1425--1443. 

\bibitem{hardt} R. Hardt and L. Simon, Nodal sets for solutions of elliptic equations. J. Diff. Geom. 30:505--522, 1989.

\bibitem{herrmann} H. Herrmann, Beitr\"age zur Theorie der Eigenwerten und Eigenfunktionen, G\"ottingen Dissertation 1932.

\bibitem{hurwitz} A. Hurwitz, \"Uber die Fourierschen Konstanten integrierbarer Funktionen, Math. Ann. 57, 425--446 (1903).
 


\bibitem{hezari} H. Hezari and Z. Wang, 
Lower bounds for volumes of nodal sets: an improvement of a result of Sogge-Zelditch. Spectral geometry, 229--235, 
Proc. Sympos. Pure Math., 84, Amer. Math. Soc., Providence, RI, 2012. 

\bibitem{hezari2} H. Hezari and C. Sogge, 
A natural lower bound for the size of nodal sets. 
Anal. PDE 5 (2012), no. 5, 1133--1137. 


\bibitem{jer} D. Jerison and G. Lebeau,
Nodal sets of sums of eigenfunctions. Harmonic analysis and partial differential equations, 223--239, 
Chicago Lectures in Math., Univ. Chicago Press, Chicago

\bibitem{jayjay} J. Jung, Bounding the number of nodal domains of eigenfunctions without singular points on the square, to appear in Israel J. Math, arXiv:1712.09457

\bibitem{lin2} F.-H. Lin, 
Nodal sets of solutions of elliptic and parabolic equations. 
Comm. Pure Appl. Math. 44 (1991), no. 3, 287--308. 

\bibitem{liouville} J. Liouville. M\'{e}moire sur le d\'{e}veloppement de fonctions ou parties de fonctions
en s\'{e}ries dont les divers termes sont assuj\'{e}tis \`a satisfaire \`a une m\^eme
\'{e}quation diff\'{e}rentielle du second ordre, contenant un param\`etre variable. J. Math. Pures Appl. 1 (1836), 253--265

\bibitem{lena} C. Lena, Courant-sharp eigenvalues of a two-dimensional torus, Comptes Rendus Mathematique 353 (2015), p. 535--539


\bibitem{logan} B. Logan, 
Integrals of high-pass functions. 
SIAM J. Math. Anal. 15 (1984), no. 2, 389--405. 


\bibitem{log2} A. Logunov, Nodal sets of Laplace eigenfunctions: polynomial upper estimates of the Hausdorff measure, Annals of Mathematics, 187,
  221-239 (2018).

\bibitem{logunov} A. Logunov, Nodal sets of Laplace eigenfunctions: proof of Nadirashvili's conjecture and of the lower bound in Yau's conjecture, Annals of Mathematics, 187,
 241--262 (2018).



\bibitem{mal} A. Logunov and E. Malinnikova, Nodal sets of Laplace eigenfunctions: estimates of the Hausdorff measure in dimension two and three, 50 Years with Hardy Spaces: A Tribute
to Victor Havin,  Operator Theory: Advances and Applications  261, Birkh\"auser.



\bibitem{mangoubi} D. Mangoubi, Local asymmetry and the inner radius of nodal domains. Comm. Part. Diff. Eqs. 33:1611--1621, 2008.




\bibitem{nadirashvili} N. Nadirashvilli, Geometry of nodal sets and multiplicity of eigenvalues, Current
Developments in Mathematics, 1997, 231--235


\bibitem{otto}  F. Otto, The geometry of dissipative evolution equations: the porous medium equation. Comm. Partial Differential Equations 26 (2001), no. 1-2, 101--174.

\bibitem{pey} R. Peyre, Non-asymptotic equivalence between $W^2$ distance and $\dot H^{-1}$, to appear in ESAIM: COCV, 2018.

\bibitem{plei} ] A. Pleijel, Remarks on Courant’s nodal line theorem, Comm. Pure Appl. Math. 9 (1956), 543--550.

\bibitem{polt}  I. Polterovich, Pleijel’s nodal domain theorem for free membranes. Proc. Amer. Math. Soc. 137 (2009), no. 3,
1021--1024.

\bibitem{pol} G. Polya, Qualitatives \"uber W\"armeausgleich, Z. Angew. Math. Mech. 13 (1933), 125--128.

\bibitem{amir} A. Sagiv and S. Steinerberger, Transport and Interface: an Uncertainty Principle for the Wasserstein distance, arXiv:1905.07450


\bibitem{sog1} C. Sogge and S. Zelditch, S. Lower bounds on the Hausdorff measure of nodal sets. Math. Res. Lett. 18:25--37. (2011)


\bibitem{sog2} C. Sogge and S. Zelditch, 
Lower bounds on the Hausdorff measure of nodal sets II. 
Math. Res. Lett. 19 (2012), no. 6, 1361--1364. 


\bibitem{stein0} S. Steinerberger, A Geometric Uncertainty Principle with an Application to Pleijel's Estimate, Annales Henri Poincare 15 (2014), 2299 -- 2319.

\bibitem{stein1} S. Steinerberger, Oscillatory functions vanish on a large set, Asian J. Math, to appear

\bibitem{stein2} S. Steinerberger, Wasserstein Distance, Fourier Series and Applications, arXiv:1803.08011

\bibitem{stein3} S. Steinerberger, Quantitative Projections in the Sturm Oscillation Theorem,  J. Math. Pures Appl., to appear.

\bibitem{sturm}  C. Sturm, M\'{e}moire sur les \'{e}quations diff\'{e}rentielles lin\'{e}aires du second ordre, J. Math. Pures Appl. 1 (1836), 106--186.

\bibitem{sturm2} C. Sturm, M\'{e}moire sur une classe d'\'{e}quations \`a diff\'{e}rences partielles, J. Math. Pures Appl. 1 (1836), 373--444.


\bibitem{villani} C. Villani, Topics in Optimal Transportation,  Graduate Studies in Mathematics, American Mathematical Society, 2003.

\bibitem{wasser} L. N. Vasershtein, Markov processes on a countable product space, describing large systems of
automata, Problemy Peredachi Informatsii, 5, 3 (1969), pp. 64--73.



\bibitem{yau}  S.-T. Yau, Problem section, Seminar on Differential Geometry, Annals of Mathematical Studies 102, Princeton, 1982, 669--706.


\bibitem{anton} A. Zettl, Sturm-Liouville Theory, Mathematical Surveys and Monographs, American Mathematical Society, 2010.
\end{thebibliography}
\end{document}